\documentclass[12pt]{article}
\newcommand{\nc}{\newcommand}  \nc{\ov}{\over} \nc{\noi}{\noindent}
\renewcommand{\sp}{\vskip1ex} \nc{\cd}{\cdots}
\nc{\ra}{\rightarrow} \nc{\iy}{\infty} 
\nc{\be}{\begin{equation}} \nc{\ee}{\end{equation}} \nc{\dl}{\delta}
\nc{\inv}{^{-1}}  \nc{\ph}{\varphi} \nc{\tl}{\widetilde}
\nc{\pht}{\tl{\ph}}
\nc{\php}{\ph_+} \nc{\phm}{\ph_-} \nc{\phtp}{\tl{\ph_+}}
\nc{\phtm}{\tl{\ph_-}}
\nc{\la}{\lambda}  \nc{\si}{\sigma}
\nc{\twotwo}[4]{\left(\begin{array}{cc}#1&#2\\&\\#3&#4\end{array}\right)}
\begin{document}

\begin{center}{\large \bf On the limit of some Toeplitz-like
determinants}\end{center}

\sp\begin{center}{{\bf Craig A.~Tracy}\\
{\it Department of Mathematics and Institute of Theoretical Dynamics\\
University of California, Davis, CA 95616, USA\\
e-mail address: tracy@itd.ucdavis.edu}}\end{center}
\begin{center}{{\bf Harold Widom}\\
{\it Department of Mathematics\\
University of California, Santa Cruz, CA 95064, USA\\
e-mail address: widom@math.ucsc.edu}}\end{center}\sp

Our Toeplitz-like matrices are of the form
\[M=(c_{p_i-q_j}),\ \ \ i,j=0,\,1,\,\cd\]
where $\{p_i\}$ and $\{q_i\}$ are sequences of integers satisfying
$p_i=q_i=i$ for $i$ sufficiently large, say for $i\ge m$. These
are a particular class of finite-rank perturbations of Toeplitz
matrices. If the $p_i$ and the $q_i$ are all different
then after rearranging the first $m$ rows and columns these become
minors of the
Toeplitz matrix $(c_{i-j})$ obtained by removing finitely many rows and
columns.

The problem is to determine the asymptotics of the determinants of large
sections of this
matrix and we do this under certain conditions.
We assume that $\{c_i\}$ is the sequence of Fourier coefficients of a
bounded
function $\ph$, so that $(c_{i-j})=T(\ph)$ in the usual notation.
We assume also that $T(\ph)$ is invertible on the space
$\ell^2(\bf{Z}^+)$, and that is
almost all. (We shall explain this below.)

For convenience in notation we consider the $(m+n)\times(m+n)$ sections
of $M$,
\[M_{m+n}=(c_{p_i-q_j}),\ \ \ i,j=0,\,1,\,\cd,\,m+n-1,\]
and denote by $T_n(\ph)$ the $n\times n$ Toeplitz matrix
$(c_{i-j})_{i,j=0,\,1,\cd,n-1}$.

If $T(\ph)$ is invertible then $\ph$ has a factorization
$\ph=\ph^-\,\ph^+$ where the
functions $(\ph^+)^{\pm1}$ and $(\ph^-)^{\pm1}$ belong to $H^2$ and
$\overline{H^2}$ respectively.\footnote{Recall that $H^2$ consists of the
$L^2$ functions whose Fourier coefficients with negative index all
vanish. The sequences of Fourier coefficients
with nonnegative indices of $\ph^+$ and $\overline{\ph^-}$ are, up to constant factors,
$T(\ph)\inv\delta$ resp. $T(\overline{\ph})\inv\delta$ where $\delta$ is the sequence
$\{1,\,0,\,0,\cd\}$.}
This is the
``Wiener-Hopf factorization'' of $\ph$. In terms of
these factors (and with subscripts denoting Fourier coefficients) our
limit formula
is\sp

\be\lim_{n\ra\iy}{\det M_{m+n}\ov\det T_n(\ph)}=
\det\,\Big(\sum_{k=1}^{\iy}
(\ph^-)_{p_i+k-m}\,(\ph^+)_{-q_j-k+m}\Big)_{i,\,j=0,\cd,m-1}.
\label{lim}\ee
(The sum in the determinant on the right side of (\ref{lim}) has only
finitely many nonzero terms since the Fourier coefficient $(\ph^-)_k$
vanishes for $k>0$
and $(\ph^+)_k$ vanishes for $k<0$.)

Under some additional hypotheses the strong Szeg\"o limit theorem gives
the asymptotics
of the Toeplitz determinant:

\[\det T_n(\ph)\sim G(\ph)^n\,E(\ph),\]
where
\[G(\ph)=\exp\left\{{1\ov2\pi}\int \log\ph(\theta)\,d\theta\right\},\ \
\
E(\ph)=\exp\left\{\sum_{k=1}^{\iy}k(\log\ph)_k\,(\log\ph)_{-k}\right\}.\]

Our result holds without these extra hypotheses. What we do need, which
is
a little stronger than
the invertibility of $T(\ph)$, is the uniform invertibility of the
$T_n(\ph)$. That
is, we require that the $T_n(\ph)$ be invertible for sufficiently large
$n$
and that the norms of the inverse matrices $T_n(\ph)\inv$ be bounded as
$n\ra\iy$.
This holds in all ``normal'' cases where $T(\ph)$ is invertible. The
known
counterexamples are not simple. (For a discussion of these points see
\cite{bs},
especially Chap. 2.) If $\ph$ has
a continuous logarithm, for example, then the $T_n(\ph)$ are uniformly
invertible.\footnote{For such a $\ph$ the Wiener-Hopf factors $\ph^-$ and $\ph^+$
are, up to constant factors, the exponentials 
of the portions of the Fourier series of $\log\ph$ corresponding to
negative resp. positive indices.}\sp

\noi{\bf Theorem}. If the $T_n(\ph)$ are uniformly invertible then
(\ref{lim}) holds.\sp

\noi{\bf Proof}. We use the fact from linear algebra that if we have a
matrix
\[M=\twotwo{A}{B}{C}{D}\]
with $A$ and $D$ square and if $D$ is invertible then
\[\det M=\det D\,\det(A-BD\inv C).\]
In our case $M=M_{m+n}$, $A=(c_{p_i-q_j})_{i,j=0,\cd,m-1}$ and $D$ is
the Toeplitz
matrix $T_n(\ph)$. If the indices for $B$
and $C$ start at 0 then their entries are given by
\be B_{i,k}=c_{p_i-m-k},\ \ \ C_{k,j}=c_{m+k-q_j},\ \ \ \
(i,j=0,\cd,m-1,\ \ k=0,\cd,n-1),
\label{BC}\ee
and we are interested in the limit
as $n\ra\iy$ of the $i,j$ entry of the $m\times m$ matrix $BT_n(\ph)\inv
C$.

Now we use the fact that if the $T_n(\ph)$ are uniformly invertible then
$T_n(\ph)\inv$
converges strongly to the infinite Toeplitz matrix
$T(\ph)\inv$. (See \cite{bs}, Prop. 2.2.) It follows that each entry of
$BT_n(\ph)\inv C$
converges to the
corresponding entry of $BT(\ph)\inv C$ where now $B$ and $C$ are the
$m\times\iy$
and $\iy\times m$ matrices, respectively, with entries given by
(\ref{BC}) but now with
$k\in\bf{Z}^+$. What remains is to show that the $i,j$ entry of
$A-BT(\ph)\inv C$ is given by the summand on the right side of
(\ref{lim}).

It is convenient to extend the range of the row or column indices of our
Toeplitz matrices
so that one or the other can run over $\bf{Z}$ rather than $\bf{Z}^+$.
So we introduce the
notations $T^{(r)}(\ph)$ and $T^{(c)}(\ph)$ for the matrices in which
the row resp. the column
index runs over $\bf{Z}$. With this notation, we see that
\[B_{i,k}=T^{(r)}(\ph)_{p_i-m,\,k}\,,\ \ \
C_{k,j}=T^{(c)}(\ph)_{k,\,q_j-m}.\]
Thus
\[(BT(\ph)\inv C)_{i,j}=(T^{(r)}(\ph)\,T(\ph)\inv
\,T^{(c)})_{p_i-m,\,q_j-m}.\]
Now it is well-known and an easy check (see \cite{bs}, Prop. 1.13) that
\[T(\ph)=T(\ph^-)\,T(\ph^+),\ \ \
T(\ph)\inv=T(\ph^+)\inv\,T(\ph^-)\inv,\]
and just as easy that
\[T^{(r)}(\ph)=T^{(r)}(\ph^-)\,T(\ph^+),\ \ \
T^{(c)}(\ph)=T(\ph^-)\,T^{(c)}(\ph^+).\]
It follows that
\[T^{(r)}(\ph)\,T(\ph)\inv\,
T^{(c)}(\ph)=T^{(r)}(\ph^-)\,T^{(c)}(\ph^+).\]
Hence
\[(BT(\ph)\inv C)_{i,j}=\sum_{k=0}^{\iy}
(\ph^-)_{p_i-m-k}\,(\ph^+)_{m+k-q_j}\]
\[=\sum_{k=-\iy}^{\iy}(\ph^-)_{p_i-m-k}\,(\ph^+)_{k-q_j+m}-
\sum_{k=-\iy}^{-1}(\ph^-)_{p_i-m-k}\,(\ph^+)_{k-q_j+m}\]
\[=\ph_{p_i-q_j}-\sum_{k=1}^{\iy}(\ph^-)_{p_i-m+k}\,(\ph^+)_{-k-q_j+m}.\]
The first term on the right is $A_{i,j}$, and so the theorem is
established.\sp\sp

\begin{center}{\bf Acknowledgments} \end{center}

Recently D. Bump and P. Diaconis \cite{bd} obtained
the asymptotics for the same determinants with an answer having a very
different form
(although it must be the same as ours) and using very different
techniques --- the
representation theory of the symmetric group. We thank them for alerting
us
to their results and sending us copies of their manuscript. We also
acknowledge support
from the National Sciance Foundation
through grants DMS-9802122 (first author) and DMS-9732687 (second
author).


\begin{thebibliography}{ss}

\bibitem{bd} D. Bump and P. Diaconis, {\it Toeplitz minors}, to appear in J. Comb. Th. A.

\bibitem{bs} A. B\"ottcher, and B. Silbermann, {\it Introduction to
large truncated
Toeplitz matrices}, Springer-Verlag, heidelberg, 1999.

\end{thebibliography}
\end{document}